\documentclass{amsart}

\usepackage{amssymb}
\usepackage{amsmath}
\usepackage{amsthm}
\usepackage{amsfonts}

\usepackage{a4wide}

\newtheorem{fundamental}{FUNDAMENTAL THEOREM}[section]
\newtheorem{theorem}{Theorem}[section]
\newtheorem{proposition}{Proposition}[section]
\newtheorem{corollary}{Corollary}[section]

\theoremstyle{definition}
\newtheorem{definition}{Definition}[section]

\newtheorem{example}{Example}[section]

\numberwithin{equation}{section}

\begin{document}

%\linenumbers

\title
[Tropen by Josef M. Hauer as orbits]{\(44\) Tropen by Josef Matthias Hauer as orbits \\ under the group of interval shifts}

\author{Helga Boyer von Berghof}
\address{Krenngasse 43\\8010 Graz\\Austria}
\email{helgaboyervonberghof@gmail.com}

\author{Daniel C. Mayer}
\address{Naglergasse 53\\8010 Graz\\Austria}
\email{quantum.algebra@icloud.com}
\urladdr{http://www.algebra.at}

%\thanks{Forschung unterst\"utzt durch den \"Osterreichischen Fonds zur F\"orderung der wissenschaftlichen Forschung (FWF): P 26008-N25}
\thanks{Research of the second author supported by the Austrian Science Fund (FWF): projects
P26008-N25, J0497-PHY, and by the European Union Research Executive Agency (EUREA): project Horizon Europe 2021--2027.}


\date{Tuesday, 17 November 2026}

\begin{abstract}
We generate a splitting of the \(924\)
monotonically ascending twelve tone half series
\(H=(H_1,\ldots,H_6)\in\lbrace 0,\ldots,11\rbrace^6\)
in \(80\) orbits \(H^G\)
under the additive cyclic group
\(G=\mathbb{Z}/12\,\mathbb{Z}\)
of \textit{transpositions} (interval shifts).
By classification according to cardinality
and forming the union of the orbits
\(H^G\cup K^G\)
of a half series \(H\) and its complement \(K\),
we obtain exactly the \(44\) \textit{Tropen}
by \textit{Josef Matthias Hauer},
completely independent of symmetry properties,
however, fundamentally due to so-called \textit{wider-gleiche} orbits.
This was not clear in previous literature.
\end{abstract}

\maketitle

%\newpage
%--------------------------------------------------------------------------------

\section{Introduction}
\label{s:Intro}

\noindent
Here, we
neither touch upon the neo-Riemannian operations (NRO)
\cite{Co1998}
nor the uniform triadic transformations (UTT)
\cite{Hk2002},
which provide a way for the description of
the chronological harmonic evolution
of a musical work of art.
However, we presuppose a familiarity
with the foundations of both of these theories.
Both are based on the concept of \textit{triads},
which we recall in section \S
\ref{s:PitchClasses}
concerning \textit{pitch classes} and
their \textit{transpositions} and inversion.
In section \S
\ref{s:HauerTropen},
we deliver an exact mathematical proof
for the precise number of \(44\) \textit{Tropen},
which J. M. Hauer
\cite{Hr1924,Hr1926}
developed as phrases
in the sense of a flexible ground material
for the composition of twelve tone plays,
independently of the chronological harmonic evolution.
This clarifies expositions in previous literature.

%\newpage
%--------------------------------------------------------------------------------

\section{Pitch classes, transpositions and inversion}
\label{s:PitchClasses}

\noindent
The \textit{scale degree} only admits
a coarse characterization of a triad,
independently of its tonic and of its mode.
The major triad \(\mathrm{XM}\) and the minor triad \(\mathrm{xm}\),
with respect to any tonic \(\mathrm{X}\), respectively \(\mathrm{x}\),
are both represented by
the triplet \(\lbrack\hat{1},\hat{3},\hat{5}\rbrack\) of scale degrees.
With sets of \textit{pitch classes}, however, we are able
to fix the tonic, for instance \(\mathrm{A}\),
and to distinguish
the major triad \(\mathrm{AM}\), \(\lbrack 0,4,7\rbrack\),
which is a symbol for \(\lbrack a,c^{\sharp},e\rbrack\),
from the minor triad \(\mathrm{am}\), \(\lbrack 0,3,7\rbrack\),
which is a symbol for \(\lbrack a,c,e\rbrack\).
An overview of scale degrees and pitch classes is given in Table
\ref{tbl:ScaleAndPitch}
with respect to the tonic \(\mathrm{A}\).

%--------------------------------------------------------------------------------

\renewcommand{\arraystretch}{1.5}

\begin{table}[ht]
\caption{Scale degrees and pitch classes}
\label{tbl:ScaleAndPitch}
\begin{center}
\begin{tabular}{|l|c|c|c|c|c|c|c|c|c|c|c|c|}
\hline
 Scale degree & \(\hat{1}\) &                & \(\hat{2}\) & \(\hat{3}\) &                & \(\hat{4}\) &                & \(\hat{5}\) & \(\hat{6}\) &                & \(\hat{7}\) &                \\
\hline
 Sharp        &             & \(a^{\sharp}\) &             &             & \(c^{\sharp}\) &             & \(d^{\sharp}\) &             &             & \(f^{\sharp}\) &             & \(g^{\sharp}\) \\
 Tonic        & \(a\)       &                & \(b\)       & \(c\)       &                & \(d\)       &                & \(e\)       & \(f\)       &                & \(g\)       &                \\
 Flat         &             & \(b^{\flat}\)  &             &             & \(d^{\flat}\)  &             & \(e^{\flat}\)  &             &             & \(g^{\flat}\)  &             & \(a^{\flat}\)  \\ 
\hline     
 Pitch class  & \(0\)       & \(1\)          & \(2\)       & \(3\)       & \(4\)          & \(5\)       & \(6\)          & \(7\)       & \(8\)       & \(9\)          & \(10\)      & \(11\)         \\      
\hline
\end{tabular}
\end{center}
\end{table}

%--------------------------------------------------------------------------------

%\noindent
As a basis for the pitch classes, we take the \textit{well tempered} scale,
which is characterized by the frequency relations 
of tones as oscillations:
\[(n:0)\quad\hat{=}\quad ((\sqrt[12]{2})^n:1), \quad \text{ for } 0\le n\le 12 \text{ (one octave)}.\]

%--------------------------------------------------------------------------------

\begin{proposition}
\label{prp:PitchClasses}
The pitch classes, taken with respect to the equally tempered scale,
together with the operation of addition modulo \(12\),
form a group \(G:=\lbrace 0,1,\ldots,11\rbrace\),
isomorphic to the cyclic group \(C_{12}\) of order \(12\),
which has the identifier \(\langle 12,2\rangle\) in the SmallGroups Database
\cite{BEO2002,BEO2005,GAP2014,MAGMA2026}.
In particular, \(G\) is an abelian (or commutative) finite group.
\end{proposition}

\begin{proof}
The neutral element of the group \(G\) is \(0\).
A generator of order \(12\) is the element \(1\).
Thus, \(G\simeq C_{12}\simeq\langle 12,2\rangle\)
is isomorphic to the cyclic group of order \(12\).
Every cyclic group is abelian.
\end{proof}

%--------------------------------------------------------------------------------

Since we would like to have a means for changing the tonic,
we introduce two kinds of permutations of \(G\), that is,
bijective mappings from \(G\) onto itself.
The first kind consists of \textit{transpositions} (interval shifts),

\begin{equation}
\label{eqn:Transpositions}
T_n:\,G\to G,\quad x\mapsto T_n(x)\equiv x+n\,(\mathrm{mod}\,12),
\end{equation}

\noindent
which move a pitch class up by \(n\) semitones, for any fixed non-negative integer \(n\).
A generic paradigm of the second kind is the \textit{inversion},

\begin{equation}
\label{eqn:Inversion}
I:\,G\to G,\quad x\mapsto I(x)\equiv -x\,(\mathrm{mod}\,12),
\end{equation}

\noindent
which generates the mirror image of a pitch class
with respect to the fixed pitch class \(0\) as the center of reflection.

%--------------------------------------------------------------------------------

\begin{theorem}
\label{thm:TranspAndInv}
The special transposition \(T_1\) by a single semitone (a small second)
and the inversion \(I\) generate a group \(\Pi:=\langle T_1,I\rangle\)
with respect to the (usual right to left) composition of mappings.
This non-abelian group \(\Pi\) is isomorphic to the dihedral group
\(D_{12}=C_{12}\rtimes C_2\) of order \(24\), that is the semidirect product
of two cyclic groups \(C_{12}\) and \(C_2\),
which has the identifier \(\langle 24,6\rangle\) in the SmallGroups Database
\cite{BEO2005}.
The transpositions alone form a normal subgroup
\(T=\lbrace T_0=1,T_1,\ldots,T_{11}\rbrace\lhd\Pi\)
isomorphic to the cyclic group \(C_{12}\simeq\langle 12,2\rangle\) of order \(12\).
The inversion generates a malnormal subgroup \(\lbrace 1,I\rbrace<\Pi\)
isomorphic to the cyclic group \(C_2\simeq\langle 2,1\rangle\) of order \(2\).
A presentation of \(\Pi\) can be given in the shape

\begin{equation}
\label{eqn:DihedralPresentation}
\Pi=\langle T_1,I\mid T_1^{12}=1,\ I^2=1,\ T_1I=IT_1^{-1}\rangle.
\end{equation}

\end{theorem}

\begin{proof}
We point out the following facts.

\begin{enumerate}

\item
The smallest power \(T_1^n\) of \(T_1\) with positive exponent \(n\ge 1\)
which coincides with the identity map is \(T_1^{12}=1\),
since \(T_1^n(x)\equiv x+n\pmod{12}\) for all \(x\in G\).

\item
The inversion \(I\) is an involution of order \(2\),
since it is its own inverse with \(I^2=1\).

\item
The inverse \(T_1^{-1}\) of \(T_1\) coincides with \(T_{11}\)
and the commutator relation can be verified by\\
\(IT_1^{-1}(x)\equiv IT_{11}(x)\equiv I(x+11)\equiv -x-11\equiv -x+1\equiv T_1(-x)\equiv T_1I(x)\pmod{12}\),\\
which is valid for any \(x\in G\).

\item
The presentation \(\Pi=\langle T_1,I\mid T_1^{12}=1,\ I^2=1,\ T_1I=IT_1^{-1}\rangle\) of the group \(\Pi\)
shows that \(\Pi\simeq D_{12}=C_{12}\rtimes C_2\) is isomorphic to the dihedral group of order \(24\).

\item
The required \textit{homomorphism}
\(\theta:\,C_2\to\mathrm{Aut}(C_{12})\),
which is tacitly underlying the semidirect product
\(D_{12}=C_{12}\rtimes C_2\),
has the trivial property \(\theta(1)=1\) and the essential property
\[\theta(I)(T_1)=I^{-1}T_1I=I^{-1}IT_1^{-1}=T_1^{-1}(=T_{11}).\qedhere\]

\end{enumerate}
\end{proof}

%--------------------------------------------------------------------------------

With transpositions alone it would be impossible to switch from a major triad to its relative minor triad.
Employing the inversion, however, we can describe such a harmonic evolution:

\begin{example}
\label{exm:TranspAndInv}
The major triad \(\mathrm{CM}\), \(\lbrack 0,4,7\rbrack\), is mapped to
\[T_7I(\lbrack 0,4,7\rbrack)=\lbrack T_7I(0),T_7I(4),T_7I(7)\rbrack=\lbrack -0+7,-4+7,-7+7\rbrack=\lbrack 7,3,0\rbrack,\]
by the permutation \(T_7I\),
that is, to the minor triad \(\mathrm{cm}\), \(\lbrack 0,3,7\rbrack\).

A drawback is that mapping the major triad \(\mathrm{DM}\), \(\lbrack 2,6,9\rbrack\),
to the minor triad \(\mathrm{dm}\), \(\lbrack 2,5,9\rbrack\),
requires a different permutation \(T_{11}I\) as the following chain of equations shows:
\[T_{11}I(\lbrack 2,6,9\rbrack)=\lbrack T_{11}I(2),T_{11}I(6),T_{11}I(9)\rbrack=\lbrack -2+11,-6+11,-9+11\rbrack=\lbrack 9,5,2\rbrack,\]
whereas our listening experience causes the expectation of a uniform triadic transformation for both mappings.

Without giving all steps, we emphasize that
\(T_7I(\lbrack 2,6,9\rbrack)=\lbrack 1,5,10\rbrack\), the \(\mathrm{b}^\flat\mathrm{m}\) triad, and
\(T_{11}I(\lbrack 0,4,7\rbrack)=\lbrack 4,7,11\rbrack\), the \(\mathrm{em}\) triad.
\end{example}

This was one of the reasons why Julian Hook developed his theory of \textit{uniform triadic transformations}
\cite{Hk2002}.
Hook received his Ph.D. in mathematics from Princeton University, NJ,
and his D.Mus. from Jacobs School of Music, Indiana University Bloomington.

David Lewin investigated contextual inversions of triads,
in particular the parallel operation PAR and the relative operation REL,
in his 1982 article
\textit{A formal theory of generalized tonal functions}
\cite{Lw1982}.
In 1987, he extended his theory in
\textit{Generalized musical intervals and transformations} (GMIT)
\cite{Lw1987}
with Riemann's leading tone exchange LTE and two kinds of shift operations,
the submediant shift MED and the subdominant shift DOM.
Brian Hyer developed Lewin's theory further
in his 1989 dissertation
\cite{Hy1989}.

%\newpage
%--------------------------------------------------------------------------------

\section{The Tropen by J. M. Hauer}
\label{s:HauerTropen}

\noindent
Now we present a precise mathematical modeling
for the deduction of the precise number of \(44\) \textit{Tropen}, which J. M. Hauer
\cite{Hr1924,Hr1926}
developed as phrases in the sense of a flexible ground material
for the composition of twelve tone plays,
independently of the chronological harmonic evolution
of a musical work of art.

%--------------------------------------------------------------------------------

\begin{definition}
\label{dfn:Reihe}
By a twelve tone series, briefly a \textit{series},
we understand an ordered dodecuplet
\(R=(R_1,\ldots,R_{12})\in Z^{12}\)
without repetitions
in the twelve tone space \(Z:=\lbrace 0,\ldots,11\rbrace\).
\end{definition}

%--------------------------------------------------------------------------------

\begin{proposition}
\label{prp:Reihe}
There are \(12!=479\,001\,600\) series in the twelve tone space \(Z\).
\end{proposition}

\begin{proof}
An ordered dodecuplet without repetitions in the twelve tone space \(Z\)
is an injective mapping of finite sets,
from \(\lbrace 1,\ldots,12\rbrace\) onto \(\lbrace 0,\ldots,11\rbrace\).
Consequently, it is also surjective and bijective,
and can be viewed as a permutation of \(Z\).
According to the counting principle, we can start by choosing
the first component \(R_1\) of the series \(R=(R_1,\ldots,R_{12})\)
in \(12\) different ways from \(Z=\lbrace 0,\ldots,11\rbrace\).
Then we continue by selecting
the second component \(R_2\) in only \(11\) ways from \(Z\setminus\lbrace R_1\rbrace\),
the third \(R_3\) in only \(10\) ways from \(Z\setminus\lbrace R_1,R_2\rbrace\),
and so on, until finally the last component \(R_{12}\) is fixed uniquely,
since it can be selected in only \(1\) way from \(Z\setminus\lbrace R_1,\ldots,R_{11}\rbrace\).
This yields a total of \(12\cdot 11\cdot 10\cdots 3\cdot 2\cdot 1=12!=479\,001\,600\) possibilities.
\end{proof}

%--------------------------------------------------------------------------------

The characteristic idea of \textbf{J. M. Hauer} was to endow the series \(R=H\dot{\cup}K\) with a structure
by splitting it disjointly into a \textbf{half series} \(H\) and its \textbf{complement} \(K=R\setminus H\).

%--------------------------------------------------------------------------------

\begin{definition}
\label{dfn:Halbreihe}
By a twelve tone half series, briefly a \textit{half series},
we understand an ordered sextet
\(H=(H_1,\ldots,H_6)\in Z^6\)
without repetitions
in the twelve tone space \(Z=\lbrace 0,\ldots,11\rbrace\).
Such a half series \(H\) is called \textit{canonical},
if it is ordered strictly monotonically ascending, that is
\(H_1<H_2<\ldots<H_6\).
\end{definition}

In the sequel, we restrict ourselves to \textit{canonical} half series \(H\),
since all other half series arise by a permutation of a canonical half series:
\(\sigma H:=(H_{\sigma(1)},\ldots,H_{\sigma(6)})\) with 
a permutation \(\sigma\) in the symmetric group \(S_6\) of degree \(6\),
which has \(6!=720\) elements.

%--------------------------------------------------------------------------------

\begin{proposition}
\label{prp:Halbreihe}
There are \(\binom{12}{6}=924\) canonical half series in the twelve tone space \(Z\).
\end{proposition}

\begin{proof}
We can choose
\(6\) objects out of \(12\) objects in
\[\binom{12}{6}=\frac{12\cdot 11\cdot 10\cdot 9\cdot 8\cdot 7}{6\cdot 5\cdot 4\cdot 3\cdot 2}
=2^2\cdot 3\cdot 7\cdot 11=924\]
ways,
without taking into account their ordering.
In particular, we can thus select 
a \textit{canonical},
that is, strictly monotonically ascending,
\textit{half series} \(H=(H_1,\ldots,H_6)\)
with \(H_1<\ldots<H_6\)
from the twelve tone space \(Z=\lbrace 0,\ldots,11\rbrace\)
in \(924\) ways.
\end{proof}

%--------------------------------------------------------------------------------

For each half series \(H\), its \textit{complement} \(K\),
which consists of \(Z\setminus\lbrace H_1,\ldots,H_6\rbrace\),
is determined uniquely and disjoint to \(H\).
Hence, there are \(924\) \textit{canonical pairs} \((H,K)\in Z^6\times Z^6\simeq Z^{12}\)
of half series in the twelve tone space \(Z\).

On each component of these pairs,
the symmetric group \(S_6\) of degree \(6\) operates by
\(\sigma H=(H_{\sigma(1)},\ldots,H_{\sigma(6)})\)
and we obtain a pair operation
\[(S_6\times S_6)\times (Z^6\times Z^6)\to (Z^6\times Z^6), ((\sigma,\tau),(H,K))\mapsto (\sigma H,\tau K).\]
Consequently, the cartesian product \(S_6\times S_6\)
of two symmetric groups \(S_6\) of degree \(6\)
operates on the pairs
by the rule \((\sigma,\tau)(H,K):=(\sigma H,\tau K)\).
Since all permutations of the twelve tone space \(Z\),
i.e., the symmetric group \(S_{12}\) of degree \(12\),
arise in this manner,
the index \((S_{12}:(S_6\times S_6))\) equals \(924\),
because \(12!=\binom{12}{6}\cdot 6!\cdot 6!\) and therefore
\[\#S_{12}=12!=479\,001\,600=
924\cdot 720\cdot 720=924\cdot 6!\cdot 6!=924\cdot\#(S_6\times S_6).\]
Thus we have proved:

\begin{corollary}
\label{cor:Index}
The cartesian product \(S_6\times S_6\)
of two symmetric groups \(S_6\) of degree \(6\)
is a subgroup
of the symmetric group \(S_{12}\) of degree \(12\)
with index
\[\bigl(S_{12}:(S_6\times S_6)\bigr)=924.\]
\end{corollary}

%--------------------------------------------------------------------------------

Although we shall devote our attention
to the complete list of Tropen by \textbf{J. M. Hauer},
written on 11 August 1948, with \textbf{numbers (No.)} in Table
\ref{tbl:Tropen},
we require yet another additional condition
for the optimal mathematical modeling
with a lexicographical order (Id.) of all orbits:

\begin{definition}
\label{dfn:Norm}
A \textit{normalization} of pairs of canonical half series
\((H,K)\in Z^6\times Z^6\) in the twelve tone space \(Z\) can be enforced
by the condition \(0\in H\), and thus \(0\not\in K\).
\end{definition}

Thus, the number of normalized canonical pairs \((H,K)\) of half series
is reduced to \(\frac{924}{2}=462\).

%--------------------------------------------------------------------------------

\begin{fundamental}
\label{thm:Fundamental}
By the operation of the additive cyclic group \(G=\mathbb{Z}/12\,\mathbb{Z}\)
of transpositions (interval shifts)
the totality of the \(924\) canonical,
that is, strictly monotonically ascending, half series
\(H=(H_1,\ldots,H_6)\) with \(H_1<\ldots<H_6\) of the twelve tone space \(Z\) 
splits into \(80\) orbits \(H^G\) with the identifiers \((\mathrm{Id}.)\) in the Table
\ref{tbl:Orbits},
and into \(44\) Tropen \(H^G\cup K^G\) as double orbits of half series \(H\) and complements \(K\)
with Hauer's numbers \((\mathbf{No.})\) in the Tables
\ref{tbl:Orbits}
and
\ref{tbl:Tropen}.
\end{fundamental}

\begin{proof}
With the aid of a program script
for the Computational Algebra System Magma
\cite{MAGMA2026}
we generate a splitting of the \(924\)
monotonically ascending \(12\)-tone half series
\(H=(H_1,\ldots,H_6)\in Z^6=\lbrace 0,\ldots,11\rbrace^6\)
into \(80\) orbits
\[H^G:=\lbrace (H_1+x,\ldots,H_6+x)\,(\mathrm{mod}\,12)\mid x\in G\rbrace\]
under the additive cyclic group
\(G=\mathbb{Z}/12\,\mathbb{Z}\)
of \textit{transpositions} (interval shifts),
respectively pitch classes.
By classification and collecting the orbits
\(H^G\cup K^G\)
of half series \(H\) and its complement \(K\)
we obtain exactly the \(44\) \textit{Tropen}
by \textit{Josef Matthias Hauer},
completely independent of symmetry properties,
hovever, essentially due to so-called \textit{wider-gleiche} orbits.
We continue the proof after the Table
\ref{tbl:Orbits}.
\end{proof}

%--------------------------------------------------------------------------------

\begin{definition}
\label{dfn:WiderGleich}
An orbit \(H^G\) of a half series \(H\) of the twelve tone space \(Z\)
under the group \(G\)
is called \textit{wider-gleich}, briefly w.g.,
if \(H^G=K^G\),
or equivalently, if \(K\in H^G\), that is,
if the complement \(K\) of the half series \(H\)
is contained in the Orbit \(H^G\) of \(H\).
\end{definition}

\newpage
%--------------------------------------------------------------------------------

\noindent
In all of the following Tables we denote
the last two elements, with two (decimal) digits, of the twelve tone space \(Z\)
by their hexadecimal symbols \(A=10\) and \(B=11\).
Then we can omit commas as separators.
In Table
\ref{tbl:Orbits},
the \(80\) orbits \(H^G\) are enumerated with identifiers (Id.)
in such a way
that the representatives \(H\) are arranged in lexicographically ascending order.
\textbf{No.} denotes Hauer's number,
and Co. the identifier of the complement \(K\) of \(H\),
\(N_1=\#H^G\), and \(N_2=\#(H^G\cup K^G)\).

\renewcommand{\arraystretch}{1.1}

\begin{table}[ht]
\caption{Identifiers of the \(80\) lexicographically ordered orbit representatives \(H\)}
\label{tbl:Orbits}
\begin{center}
\begin{tabular}{|r|c|c|r|r|c||r|c|c|r|r|c|}
\hline
    Id. & \(H\)      & \(K\)      &    Co. &    \textbf{No.} & \(N_1,N_2\) &     Id. & \(H\)      & \(K\)      &    Co. &    \textbf{No.} & \(N_1,N_2\) \\
\hline
  \(1\) & \(012345\) & \(6789AB\) &  \(1\) &  \(\mathbf{1}\) & \(12,12\)   &  \(41\) & \(012578\) & \(3469AB\) & \(48\) & \(\mathbf{30}\) & \(12,24\)   \\
  \(2\) & \(012346\) & \(5789AB\) &  \(6\) &  \(\mathbf{2}\) & \(12,24\)   &  \(42\) & \(012579\) & \(3468AB\) & \(67\) & \(\mathbf{36}\) & \(12,24\)   \\
  \(3\) & \(012347\) & \(5689AB\) & \(21\) &  \(\mathbf{5}\) & \(12,24\)   &  \(43\) & \(01257A\) & \(34689B\) & \(64\) & \(\mathbf{38}\) & \(12,24\)   \\
  \(4\) & \(012348\) & \(5679AB\) & \(22\) &  \(\mathbf{7}\) & \(12,24\)   &  \(44\) & \(012589\) & \(3467AB\) & \(59\) & \(\mathbf{15}\) & \(12,24\)   \\
  \(5\) & \(012349\) & \(5678AB\) &  \(7\) &  \(\mathbf{6}\) & \(12,24\)   &  \(45\) & \(01258A\) & \(34679B\) & \(56\) & \(\mathbf{33}\) & \(12,24\)   \\
  \(6\) & \(01234A\) & \(56789B\) &  \(2\) &  \(\mathbf{2}\) & \(12,24\)   &  \(46\) & \(01259A\) & \(34679B\) & \(46\) & \(\mathbf{28}\) & \(12,12\)   \\  
  \(7\) & \(012356\) & \(4789AB\) &  \(5\) &  \(\mathbf{6}\) & \(12,24\)   &  \(47\) & \(012678\) & \(3459AB\) & \(47\) &  \(\mathbf{4}\) & \(6,6\)     \\
  \(8\) & \(012357\) & \(4689AB\) & \(20\) &  \(\mathbf{9}\) & \(12,24\)   &  \(48\) & \(012679\) & \(3458AB\) & \(41\) & \(\mathbf{30}\) & \(12,24\)   \\
  \(9\) & \(012358\) & \(4679AB\) & \(54\) & \(\mathbf{23}\) & \(12,24\)   &  \(49\) & \(01267A\) & \(34589B\) & \(38\) & \(\mathbf{24}\) & \(12,24\)   \\
 \(10\) & \(012359\) & \(4678AB\) & \(26\) & \(\mathbf{19}\) & \(12,24\)   &  \(50\) & \(012689\) & \(3457AB\) & \(31\) & \(\mathbf{25}\) & \(12,24\)   \\
 \(11\) & \(01235A\) & \(46789B\) & \(11\) & \(\mathbf{10}\) & \(12,12\)   &  \(51\) & \(01268A\) & \(34579B\) & \(28\) & \(\mathbf{11}\) & \(12,24\)   \\ 
 \(12\) & \(012367\) & \(4589AB\) & \(19\) &  \(\mathbf{3}\) & \(12,24\)   &  \(52\) & \(01269A\) & \(34578B\) & \(24\) & \(\mathbf{26}\) & \(12,24\)   \\
 \(13\) & \(012368\) & \(4579AB\) & \(53\) & \(\mathbf{18}\) & \(12,24\)   &  \(53\) & \(01278A\) & \(34569B\) & \(13\) & \(\mathbf{18}\) & \(12,24\)   \\
 \(14\) & \(012369\) & \(4578AB\) & \(55\) & \(\mathbf{14}\) & \(12,24\)   &  \(54\) & \(01279A\) & \(34568B\) &  \(9\) & \(\mathbf{23}\) & \(12,24\)   \\
 \(15\) & \(01236A\) & \(45789B\) & \(36\) & \(\mathbf{21}\) & \(12,24\)   &  \(55\) & \(013467\) & \(2589AB\) & \(14\) & \(\mathbf{14}\) & \(12,24\)   \\ 
 \(16\) & \(012378\) & \(4569AB\) & \(37\) &  \(\mathbf{8}\) & \(12,24\)   &  \(56\) & \(013468\) & \(2579AB\) & \(45\) & \(\mathbf{33}\) & \(12,24\)   \\
 \(17\) & \(012379\) & \(4568AB\) & \(27\) & \(\mathbf{22}\) & \(12,24\)   &  \(57\) & \(013469\) & \(2578AB\) & \(61\) & \(\mathbf{34}\) & \(12,24\)   \\
 \(18\) & \(01237A\) & \(45689B\) & \(23\) & \(\textbf{20}\) & \(12,24\)   &  \(58\) & \(01346A\) & \(25789B\) & \(33\) & \(\textbf{32}\) & \(12,24\)   \\ 
 \(19\) & \(012389\) & \(4567AB\) & \(12\) &  \(\textbf{3}\) & \(12,24\)   &  \(59\) & \(013478\) & \(2569AB\) & \(44\) & \(\textbf{15}\) & \(12,24\)   \\
 \(20\) & \(01238A\) & \(45679B\) &  \(8\) &  \(\textbf{9}\) & \(12,24\)   &  \(60\) & \(013479\) & \(2568AB\) & \(65\) & \(\textbf{35}\) & \(12,24\)   \\
 \(21\) & \(01239A\) & \(45678B\) &  \(3\) &  \(\textbf{5}\) & \(12,24\)   &  \(61\) & \(01347A\) & \(25689B\) & \(57\) & \(\textbf{34}\) & \(12,24\)   \\
 \(22\) & \(012456\) & \(3789AB\) &  \(4\) &  \(\textbf{7}\) & \(12,24\)   &  \(62\) & \(013489\) & \(2567AB\) & \(39\) & \(\textbf{16}\) & \(12,24\)   \\
 \(23\) & \(012457\) & \(3689AB\) & \(18\) & \(\textbf{20}\) & \(12,24\)   &  \(63\) & \(01348A\) & \(25679B\) & \(29\) & \(\textbf{31}\) & \(12,24\)   \\
 \(24\) & \(012458\) & \(3679AB\) & \(52\) & \(\textbf{26}\) & \(12,24\)   &  \(64\) & \(013568\) & \(2479AB\) & \(43\) & \(\textbf{38}\) & \(12,24\)   \\
 \(25\) & \(012459\) & \(3678AB\) & \(25\) & \(\textbf{29}\) & \(12,12\)   &  \(65\) & \(013569\) & \(2478AB\) & \(60\) & \(\textbf{35}\) & \(12,24\)   \\
 \(26\) & \(01245A\) & \(36789B\) & \(10\) & \(\textbf{19}\) & \(12,24\)   &  \(66\) & \(01356A\) & \(24789B\) & \(32\) & \(\textbf{37}\) & \(12,24\)   \\
 \(27\) & \(012467\) & \(3589AB\) & \(17\) & \(\textbf{22}\) & \(12,24\)   &  \(67\) & \(013578\) & \(2469AB\) & \(42\) & \(\textbf{36}\) & \(12,24\)   \\ 
 \(28\) & \(012468\) & \(3579AB\) & \(51\) & \(\textbf{11}\) & \(12,24\)   &  \(68\) & \(013579\) & \(2468AB\) & \(79\) & \(\textbf{43}\) & \(12,24\)   \\
 \(29\) & \(012469\) & \(3578AB\) & \(63\) & \(\textbf{31}\) & \(12,24\)   &  \(69\) & \(01357A\) & \(24689B\) & \(75\) & \(\textbf{42}\) & \(12,24\)   \\
 \(30\) & \(01246A\) & \(35789B\) & \(35\) & \(\textbf{12}\) & \(12,24\)   &  \(70\) & \(013589\) & \(2467AB\) & \(77\) & \(\textbf{13}\) & \(12,24\)   \\
 \(31\) & \(012478\) & \(3569AB\) & \(50\) & \(\textbf{25}\) & \(12,24\)   &  \(71\) & \(01358A\) & \(24679B\) & \(71\) & \(\textbf{41}\) & \(12,12\)   \\ 
 \(32\) & \(012479\) & \(3568AB\) & \(66\) & \(\textbf{37}\) & \(12,24\)   &  \(72\) & \(013679\) & \(2458AB\) & \(78\) & \(\textbf{39}\) & \(6,12\)    \\
 \(33\) & \(01247A\) & \(35689B\) & \(58\) & \(\textbf{32}\) & \(12,24\)   &  \(73\) & \(01367A\) & \(24589B\) & \(74\) & \(\textbf{40}\) & \(12,24\)   \\
 \(34\) & \(012489\) & \(3567AB\) & \(40\) & \(\textbf{27}\) & \(12,24\)   &  \(74\) & \(013689\) & \(2457AB\) & \(73\) & \(\textbf{40}\) & \(12,24\)   \\ 
 \(35\) & \(01248A\) & \(35679B\) & \(30\) & \(\textbf{12}\) & \(12,24\)   &  \(75\) & \(01368A\) & \(24579B\) & \(69\) & \(\textbf{42}\) & \(12,24\)   \\
 \(36\) & \(01249A\) & \(35678B\) & \(15\) & \(\textbf{21}\) & \(12,24\)   &  \(76\) & \(014589\) & \(2367AB\) & \(76\) & \(\textbf{17}\) & \(4,4\)     \\
 \(37\) & \(012567\) & \(3489AB\) & \(16\) &  \(\textbf{8}\) & \(12,24\)   &  \(77\) & \(01458A\) & \(23679B\) & \(70\) & \(\textbf{13}\) & \(12,24\)   \\
 \(38\) & \(012568\) & \(3479AB\) & \(49\) & \(\textbf{24}\) & \(12,24\)   &  \(78\) & \(01467A\) & \(23589B\) & \(72\) & \(\textbf{39}\) & \(6,12\)    \\
 \(39\) & \(012569\) & \(3478AB\) & \(62\) & \(\textbf{16}\) & \(12,24\)   &  \(79\) & \(01468A\) & \(23579B\) & \(68\) & \(\textbf{43}\) & \(12,24\)   \\
 \(40\) & \(01256A\) & \(34789B\) & \(34\) & \(\textbf{27}\) & \(12,24\)   &  \(80\) & \(02468A\) & \(13579B\) & \(80\) & \(\textbf{44}\) & \(2,2\)     \\
\hline
\end{tabular}
\end{center}
\end{table}

\newpage
%--------------------------------------------------------------------------------

\noindent
\textbf{Continuation of the proof.}
The additive cyclic group \(G=(\mathbb{Z}/12\,\mathbb{Z},+)\)
of \textit{interval shifts} (transpositions)
operates on the canonical half series \(H\)
by addition modulo \(12\) and subsequent strictly monotonically ascending ordering.
Hereby their totality splits in
\(80\) orbits \(H^G\) with various cardinalities.
Most of these orbits, namely \(75\) among them, have \(12\) elements,
that is, trivial \textit{stabilizer}
\(S_H:=\lbrace g\in G\mid H^g=H\rbrace\).
The orbits of five distinguished half series possess smaller numbers of elements,
as we can see in the Table
\ref{tbl:KleineOrbits}.

%--------------------------------------------------------------------------------

\renewcommand{\arraystretch}{1.5}

\begin{table}[ht]
\caption{Five distinguished small orbits \(H^G\) of half series \(H\)}
\label{tbl:KleineOrbits}
\begin{center}
\begin{tabular}{|c|c|c|c|c||r|c|}
\hline
 Orbit Id. & \(H\)      & \(K\)      & \(\#H^G\) & \(\#S_H\) & Hauer \textbf{No.} & w.g.       \\
\hline
 \(80\)    & \(02468A\) & \(13579B\) & \(2\)     & \(6\)     & \(\mathbf{44}\)    & \(\times\) \\
 \(76\)    & \(014589\) & \(2367AB\) & \(4\)     & \(3\)     & \(\mathbf{17}\)    & \(\times\) \\
 \(47\)    & \(012678\) & \(3459AB\) & \(6\)     & \(2\)     & \(\mathbf{4}\)     & \(\times\) \\      
 \(72\)    & \(013679\) & \(2458AB\) & \(6\)     & \(2\)     & \(\mathbf{39}\)    &            \\ 
 \(78\)    & \(01467A\) & \(23589B\) & \(6\)     & \(2\)     & \(\mathbf{39}\)    &            \\ 
\hline
 \multicolumn{3}{|c|}{all the other \(75\) orbits} & \(12\) & \(1\) &             &            \\
\hline
\end{tabular}
\end{center}
\end{table}

%--------------------------------------------------------------------------------

\noindent
Wider-gleiche (w.g.) orbits are designated by the symbol \(\times\).
In Table
\ref{tbl:Orbits},
w.g. are those with the property \(N_1=N_2\).
But there are also \(5\) wider-gleiche among the \(75\) orbits with maximal cardinality \(12\),
how Table
\ref{tbl:Zensus}
shows.
Herewith, we are now in the position to present
the complete census of all Tropen as double orbits \(H^G\cup K^G\)
of canonical half series \(H\) and their complements \(K\)
under the group \(G\) of transpositions (interval shifts).

\renewcommand{\arraystretch}{1.5}

\begin{table}[ht]
\caption{Census of Tropen as double orbits \(H^G\cup K^G\)}
\label{tbl:Zensus}
\begin{center}
\begin{tabular}{|r|r||c|c|r|r|c|c||c|}
\hline
        Id. & \textbf{No.}    & \(H\)      & \(K\)      & \(\#H^G\) & \(\#(H^G\cup K^G)\) & w.g.       & \(\#S_H\) & Z\"ahlung \\
\hline
     \(80\) & \(\mathbf{44}\) & \(02468A\) & \(13579B\) &     \(2\) &               \(2\) & \(\times\) & \(6\)     & \(1\)     \\
     \(76\) & \(\mathbf{17}\) & \(014589\) & \(2367AB\) &     \(4\) &               \(4\) & \(\times\) & \(3\)     & \(1\)     \\
\hline
     \(47\) &  \(\mathbf{4}\) & \(012678\) & \(3459AB\) &     \(6\) &               \(6\) & \(\times\) & \(2\)     & \(1\)     \\      
 \(72\&78\) & \(\mathbf{39}\) & \(013679\) & \(2458AB\) &     \(6\) &              \(12\) &            & \(2\)     & \(2/1\)   \\ 
\hline
      \(1\) &  \(\mathbf{1}\) & \(012345\) & \(6789AB\) &    \(12\) &              \(12\) & \(\times\) & \(1\)     & \(1\)     \\
     \(11\) & \(\mathbf{10}\) & \(023457\) & \(1689AB\) &    \(12\) &              \(12\) & \(\times\) & \(1\)     & \(1\)     \\
     \(25\) & \(\mathbf{29}\) & \(0237AB\) & \(145689\) &    \(12\) &              \(12\) & \(\times\) & \(1\)     & \(1\)     \\
     \(46\) & \(\mathbf{28}\) & \(23567A\) & \(01489B\) &    \(12\) &              \(12\) & \(\times\) & \(1\)     & \(1\)     \\
     \(71\) & \(\mathbf{41}\) & \(02357A\) & \(14689B\) &    \(12\) &              \(12\) & \(\times\) & \(1\)     & \(1\)     \\
\hline
 \multicolumn{4}{|c|}{all the other \(35\) Tropen}      &    \(12\) &              \(24\) &            & \(1\)     & \(70/35\) \\
\hline
 \multicolumn{8}{|c||}{Total}                                                                                      & \(80/44\) \\
\hline
\end{tabular}
\end{center}
\end{table}

\noindent
Thus, there are
\(3\) Tropen of wider-gleiche small Orbits,
\(1\) Trope of two further small orbits,
\(5\) Tropen of wider-gleiche big orbits,
and the dominating contribution by
\(35\) Tropen of the \(70\) further big orbits,
in the total sum:
\(3+1+5+35=44\) Tropen from \(3+2+5+70=80\) orbits. QED.

\newpage
%--------------------------------------------------------------------------------

\noindent
In Table
\ref{tbl:Tropen},
the \(44\) Tropen
\((H^G\cup K^G)\,(\mathrm{mod}\,G)\)
are listed
with fixed numbering \textbf{(No.)}
written by \textbf{Josef Matthias Hauer}
on Wednesday, 11 August 1948,
and put in correspondence with canonical half series
by selection of
representing half series \(H\) and complements \(K\),
and with individual identifiers (w.g. orbits) or pairs of identifiers (Id.)

\renewcommand{\arraystretch}{1.5}

\begin{table}[ht]
\caption{The \(44\) Tropen by Josef M. Hauer of 11 August 1948}
\label{tbl:Tropen}
\begin{center}
\begin{tabular}{|r|c|c|c|c||r|c|c|c|c|}
\hline
 \textbf{No.}    & \(H\)      & \(K\)      &   Id. & w.g.       & \textbf{No.}    & \(H\)      & \(K\)       &   Id. & w.g.       \\
\hline
  \(\mathbf{1}\) & \(012345\) & \(6789AB\) &     1 & \(\times\) & \(\mathbf{23}\) & \(023567\) & \(1489AB\)  &  9,54 &            \\
  \(\mathbf{2}\) & \(012346\) & \(5789AB\) &   2,6 &            & \(\mathbf{24}\) & \(01267A\) & \(34589B\)  & 38,49 &            \\
  \(\mathbf{3}\) & \(012367\) & \(4589AB\) & 12,19 &            & \(\mathbf{25}\) & \(023678\) & \(1459AB\)  & 31,50 &            \\      
  \(\mathbf{4}\) & \(012678\) & \(3459AB\) &    47 & \(\times\) & \(\mathbf{26}\) & \(012458\) & \(3679AB\)  & 24,52 &            \\ 
  \(\mathbf{5}\) & \(012347\) & \(5689AB\) &  3,21 &            & \(\mathbf{27}\) & \(04589A\) & \(12367B\)  & 34,40 &            \\ 
  \(\mathbf{6}\) & \(012356\) & \(4789AB\) &   5,7 &            & \(\mathbf{28}\) & \(23567A\) & \(01489B\)  &    46 & \(\times\) \\ 
  \(\mathbf{7}\) & \(012456\) & \(3789AB\) &  4,22 &            & \(\mathbf{29}\) & \(0237AB\) & \(145689\)  &    25 & \(\times\) \\ 
  \(\mathbf{8}\) & \(012567\) & \(3489AB\) & 16,37 &            & \(\mathbf{30}\) & \(012578\) & \(3469AB\)  & 41,48 &            \\
  \(\mathbf{9}\) & \(012357\) & \(4689AB\) &  8,20 &            & \(\mathbf{31}\) & \(02356A\) & \(14789B\)  & 29,63 &            \\
 \(\mathbf{10}\) & \(023457\) & \(1689AB\) &    11 & \(\times\) & \(\mathbf{32}\) & \(01346A\) & \(25789B\)  & 33,58 &            \\
 \(\mathbf{11}\) & \(012468\) & \(3579AB\) & 28,51 &            & \(\mathbf{33}\) & \(013468\) & \(2579AB\)  & 45,56 &            \\
 \(\mathbf{12}\) & \(023468\) & \(1579AB\) & 30,35 &            & \(\mathbf{34}\) & \(013469\) & \(2578AB\)  & 57,61 &            \\
 \(\mathbf{13}\) & \(02367A\) & \(14589B\) & 70,77 &            & \(\mathbf{35}\) & \(013479\) & \(2568AB\)  & 60,65 &            \\      
 \(\mathbf{14}\) & \(013467\) & \(2589AB\) & 14,55 &            & \(\mathbf{36}\) & \(01568A\) & \(23479B\)  & 42,67 &            \\ 
 \(\mathbf{15}\) & \(013478\) & \(2569AB\) & 44,59 &            & \(\mathbf{37}\) & \(023578\) & \(1469AB\)  & 32,66 &            \\ 
 \(\mathbf{16}\) & \(014578\) & \(2369AB\) & 39,62 &            & \(\mathbf{38}\) & \(013568\) & \(2479AB\)  & 43,64 &            \\ 
 \(\mathbf{17}\) & \(014589\) & \(2367AB\) &    76 & \(\times\) & \(\mathbf{39}\) & \(013679\) & \(2458AB\)  & 72,78 &            \\ 
 \(\mathbf{18}\) & \(01278A\) & \(34569B\) & 13,53 &            & \(\mathbf{40}\) & \(014679\) & \(2358AB\)  & 73,74 &            \\
 \(\mathbf{19}\) & \(01245A\) & \(36789B\) & 10,26 &            & \(\mathbf{41}\) & \(02357A\) & \(14689B\)  &    71 & \(\times\) \\
 \(\mathbf{20}\) & \(012457\) & \(3689AB\) & 18,23 &            & \(\mathbf{42}\) & \(023579\) & \(1468AB\)  & 69,75 &            \\
 \(\mathbf{21}\) & \(05689A\) & \(12347B\) & 15,36 &            & \(\mathbf{43}\) & \(013579\) & \(2468AB\)  & 68,79 &            \\
 \(\mathbf{22}\) & \(02389A\) & \(14567B\) & 17,27 &            & \(\mathbf{44}\) & \(02468A\) & \(13579B\)  &    80 & \(\times\) \\
\hline
\end{tabular}
\end{center}
\end{table}

%\newpage
%--------------------------------------------------------------------------------

\section{New lexicographic arrangement of the \(44\) Tropen}
\label{s:LexTropen}

\noindent
In Table
\ref{tbl:LexTropen},
the \(44\) Tropen
\((H^G\cup K^G)\,(\mathrm{mod}\,G)\)
with a \textbf{new enumeration} \textbf{(Tr.)} are arranged in such a manner
that the representing half series \(H\) are
\textit{lexicographically ascending} and always \textit{normalized} by \(0\in H\).
Instead of the complement \(K\) of \(H\), which can be found in the Table
\ref{tbl:Orbits},
we always give the normalized representative \(H(K)\) of the orbit \(K^G\).
With pairs of orbit identifiers (Id., Co.)
and the Hauer Tropen number \textbf{(No.)}
the correspondence to the earlier Tables is provided.
The cardinalities \(N_1,N_2\) are taken from the Table
\ref{tbl:Orbits}. 

\renewcommand{\arraystretch}{1.5}

\begin{table}[ht]
\caption{New identifiers Tr. of \(44\) Tropen (representatives \(H\) arranged lexicographically)}
\label{tbl:LexTropen}
\begin{center}
\begin{tabular}{|r|r|c|c|r|r|c||r|r|c|c|r|r|c|}
\hline
    \textbf{Tr.} &    Id. & \(H\)      & \(H(K)\)   &    Co. &    \textbf{No.} & \(N_1,N_2\) &    \textbf{Tr.} &    Id. & \(H\)      & \(H(K)\)   &    Co. &    \textbf{No.} & \(N_1,N_2\) \\
\hline
  \(\mathbf{1}\) &  \(1\) & \(012345\) & \(012345\) &  \(1\) &  \(\mathbf{1}\) & \(12,12\)   & \(\mathbf{23}\) & \(32\) & \(012479\) & \(01356A\) & \(66\) & \(\mathbf{37}\) & \(12,24\)   \\
  \(\mathbf{2}\) &  \(2\) & \(012346\) & \(01234A\) &  \(6\) &  \(\mathbf{2}\) & \(12,24\)   & \(\mathbf{24}\) & \(33\) & \(01247A\) & \(01346A\) & \(58\) & \(\mathbf{32}\) & \(12,24\)   \\
  \(\mathbf{3}\) &  \(3\) & \(012347\) & \(01239A\) & \(21\) &  \(\mathbf{5}\) & \(12,24\)   & \(\mathbf{25}\) & \(34\) & \(012489\) & \(01256A\) & \(40\) & \(\mathbf{27}\) & \(12,24\)   \\
  \(\mathbf{4}\) &  \(4\) & \(012348\) & \(012456\) & \(22\) &  \(\mathbf{7}\) & \(12,24\)   & \(\mathbf{26}\) & \(38\) & \(012568\) & \(01267A\) & \(49\) & \(\mathbf{24}\) & \(12,24\)   \\
  \(\mathbf{5}\) &  \(5\) & \(012349\) & \(012356\) &  \(7\) &  \(\mathbf{6}\) & \(12,24\)   & \(\mathbf{27}\) & \(39\) & \(012569\) & \(013489\) & \(62\) & \(\mathbf{16}\) & \(12,24\)   \\
  \(\mathbf{6}\) &  \(8\) & \(012357\) & \(01238A\) & \(20\) &  \(\mathbf{9}\) & \(12,24\)   & \(\mathbf{28}\) & \(41\) & \(012578\) & \(012679\) & \(48\) & \(\mathbf{30}\) & \(12,24\)   \\  
  \(\mathbf{7}\) &  \(9\) & \(012358\) & \(01279A\) & \(54\) & \(\mathbf{23}\) & \(12,24\)   & \(\mathbf{29}\) & \(42\) & \(012579\) & \(013578\) & \(67\) & \(\mathbf{36}\) & \(12,24\)   \\
  \(\mathbf{8}\) & \(10\) & \(012359\) & \(01245A\) & \(26\) & \(\mathbf{19}\) & \(12,24\)   & \(\mathbf{30}\) & \(43\) & \(01257A\) & \(013568\) & \(64\) & \(\mathbf{38}\) & \(12,24\)   \\
  \(\mathbf{9}\) & \(11\) & \(01235A\) & \(01235A\) & \(11\) & \(\mathbf{10}\) & \(12,12\)   & \(\mathbf{31}\) & \(44\) & \(012589\) & \(013478\) & \(59\) & \(\mathbf{15}\) & \(12,24\)   \\
 \(\mathbf{10}\) & \(12\) & \(012367\) & \(012389\) & \(19\) &  \(\mathbf{3}\) & \(12,24\)   & \(\mathbf{32}\) & \(45\) & \(01258A\) & \(013468\) & \(56\) & \(\mathbf{33}\) & \(12,24\)   \\
 \(\mathbf{11}\) & \(13\) & \(012368\) & \(01278A\) & \(53\) & \(\mathbf{18}\) & \(12,24\)   & \(\mathbf{33}\) & \(46\) & \(01259A\) & \(01259A\) & \(46\) & \(\mathbf{28}\) & \(12,12\)   \\ 
 \(\mathbf{12}\) & \(14\) & \(012369\) & \(013467\) & \(55\) & \(\mathbf{14}\) & \(12,24\)   & \(\mathbf{34}\) & \(47\) & \(012678\) & \(012678\) & \(47\) &  \(\mathbf{4}\) & \(6,6\)     \\
 \(\mathbf{13}\) & \(15\) & \(01236A\) & \(01249A\) & \(36\) & \(\mathbf{21}\) & \(12,24\)   & \(\mathbf{35}\) & \(57\) & \(013469\) & \(01347A\) & \(61\) & \(\mathbf{34}\) & \(12,24\)   \\ 
 \(\mathbf{14}\) & \(16\) & \(012378\) & \(012567\) & \(37\) &  \(\mathbf{8}\) & \(12,24\)   & \(\mathbf{36}\) & \(60\) & \(013479\) & \(013569\) & \(65\) & \(\mathbf{35}\) & \(12,24\)   \\
 \(\mathbf{15}\) & \(17\) & \(012379\) & \(012467\) & \(27\) & \(\mathbf{22}\) & \(12,24\)   & \(\mathbf{37}\) & \(68\) & \(013579\) & \(01468A\) & \(79\) & \(\mathbf{43}\) & \(12,24\)   \\
 \(\mathbf{16}\) & \(18\) & \(01237A\) & \(012457\) & \(23\) & \(\mathbf{20}\) & \(12,24\)   & \(\mathbf{38}\) & \(69\) & \(01357A\) & \(01368A\) & \(75\) & \(\mathbf{42}\) & \(12,24\)   \\ 
 \(\mathbf{17}\) & \(24\) & \(012458\) & \(01269A\) & \(52\) & \(\mathbf{26}\) & \(12,24\)   & \(\mathbf{39}\) & \(70\) & \(013589\) & \(01458A\) & \(77\) & \(\mathbf{13}\) & \(12,24\)   \\
 \(\mathbf{18}\) & \(25\) & \(012459\) & \(012459\) & \(25\) & \(\mathbf{29}\) & \(12,12\)   & \(\mathbf{40}\) & \(71\) & \(01358A\) & \(01358A\) & \(71\) & \(\mathbf{41}\) & \(12,12\)   \\ 
 \(\mathbf{19}\) & \(28\) & \(012468\) & \(01268A\) & \(51\) & \(\mathbf{11}\) & \(12,24\)   & \(\mathbf{41}\) & \(72\) & \(013679\) & \(01467A\) & \(78\) & \(\mathbf{39}\) & \(6,12\)    \\
 \(\mathbf{20}\) & \(29\) & \(012469\) & \(01348A\) & \(63\) & \(\mathbf{31}\) & \(12,24\)   & \(\mathbf{42}\) & \(73\) & \(01367A\) & \(013689\) & \(74\) & \(\mathbf{40}\) & \(12,24\)   \\
 \(\mathbf{21}\) & \(30\) & \(01246A\) & \(01248A\) & \(35\) & \(\mathbf{12}\) & \(12,24\)   & \(\mathbf{43}\) & \(76\) & \(014589\) & \(014589\) & \(76\) & \(\mathbf{17}\) & \(4,4\)     \\
 \(\mathbf{22}\) & \(31\) & \(012478\) & \(012689\) & \(50\) & \(\mathbf{25}\) & \(12,24\)   & \(\mathbf{44}\) & \(80\) & \(02468A\) & \(02468A\) & \(80\) & \(\mathbf{44}\) & \(2,2\)     \\
\hline
\end{tabular}
\end{center}
\end{table}

%\newpage
%--------------------------------------------------------------------------------

\section{An interesting application}
\label{s:Tristan}

\noindent
As an example, we determine the Trope
which is defined by the chord in bar 1 of J. M. Hauer's
composition \textbf{Nomos, op. 1,} for piano solo.
This chord is said to arise in Richard Wagner's opera, respectively prelude,
\textbf{Tristan and Isolde}.
It is therefore called the \textbf{Tristan chord},
\begin{center}
\(\mathrm{D}_1\) -- \(\mathrm{B}^\flat_1\) -- \(\mathrm{F}\) -- \(C^\flat\)-- \(e^\flat\) -- \(a^\flat\).
\end{center}
However, from our perspective it is a superposition of
a \(\mathrm{B}^\flat\) Major chord and an \(a^\flat\) minor chord,
which is characteristic for \textit{spanish flamenco}, due to \textit{mauretan gipsies}.
With respect to the tonic \(\mathrm{A}\),
the Tristan chord is the half series \((51826B)\),
respectively the canonical half series \((12568B)\).
In order to find the corresponding normalized representative
either as half series \(H\) or as complement \(H(K)\) in Table
\ref{tbl:LexTropen},
we apply a transposition by \(-5\hat{=}+7\)
(a shift down by a quart interval, or up by a quint interval).
So, with respect to the tonic \(\mathrm{D}\),
the Tristan chord is then the half series \((890136)\),
respectively the canonical half series \((013689)\).
This is the complement \(H(K)\), Co. 74, of the half series \(H=(01367A)\), Id. 73.
Consequently, the \textbf{Tristan chord} is precisely Hauer's Trope \textbf{No. 40},
respectively \textbf{Tr. 42} in new lexicographical arrangement.

%\newpage
%--------------------------------------------------------------------------------

\section{Acknowledgements}
\label{s:Thanks}

The authors would like to express their gratitude to Herrn
\textit{em. o. Univ.-Prof. Gerhard Zeller}
from the NOMOS Bildungs- und Forschungs-Institut
f\"ur ganzheitlich-harmonikale Strukturformen
Zahl-Farbe-Klang im Werk von Hans Florey,
and Frau 
\textit{Dr. Erica Bisesi}
from the Center of Systematic Music Theory,
both at the University of Music and Arts in Graz, Austria,
for stimulating our interest in studying the
connections between musical expression
and mathematical structures.

The second author gratefully acknowledges that his research is supported by
the Austrian Science Fund (FWF): projects J0497-PHY, P26008-N25,
and by the Research Executive Agency of the European Union (EUREA):
project Horizon Europe 2021--2027.

%\newpage
%--------------------------------------------------------------------------------

\section{Data availability}
\label{s:Data}

\noindent
Magma program scripts and output
may be requested from the second author via email.

%\newpage
%--------------------------------------------------------------------------------

\section{Conflict of interests}
\label{s:Conflict}

\noindent
The authors declare that no conflict of interests arises from this article.

%\newpage
%--------------------------------------------------------------------------------

%--------------------------------------------------------------------------------


\begin{thebibliography}{XX}
%
\bibitem{BEO2002}
%\bysame
H. U. Besche, B. Eick, and E. A. O'Brien,
\textit{A millennium project: constructing small groups},
Int. J. Algebra Comput.
\textbf{12}
(2002),
623--644.
%
\bibitem{BEO2005}
%\bysame
H. U. Besche, B. Eick, and E. A. O'Brien,
\textit{The SmallGroups Library --- a Library of Groups of Small Order},
2005,
an accepted and refereed GAP 4 package, available also in MAGMA.
%
\bibitem{Co1998}
R. Cohn,
\textit{Introduction to neo-Riemannian theory:
A survey and a historical perspective},
Journal of Music Theory
\textbf{42}
(1998),
no. 2,
167--180.
%
\bibitem{GAP2014}
GAP Developer Group,
\textit{GAP -- Groups, Algorithms, and Programming -- a System for Computational Discrete Algebra},
Version 4.16.1,
Aachen, Braunschweig, Fort Collins, St. Andrews,
2026, \\
\texttt{http://www.gap-system.org}.
%
\bibitem{Hr1924}
J. M. Hauer,
\textit{Die Tropen},
pp. 18--21 in \lq\lq Musikbl\"atter des Anbruch\rq\rq,
Universal Edition, Wien, 1924.
%
\bibitem{Hr1926}
J. M. Hauer,
\textit{Zw\"olftontechnik. Die Lehre von den Tropen},
Universal Edition, Wien, 1926.
%
\bibitem{Hk2002}
J. Hook,
\textit{Uniform triadic transformations},
Journal of Music Theory
\textbf{46}
(2002),
no. 1--2,
57--126.
%
\bibitem{Hy1989}
B. Hyer,
\textit{Tonal intuitions in \lq Tristan und Isolde\rq},
Ph.D. thesis,
Yale University,
New Haven,
CT,
1989.
%
\bibitem{Lw1982}
D. Lewin,
\textit{A formal theory of generalized tonal functions},
Journal of Music Theory
\textbf{26}
(1982),
23--60.
%
\bibitem{Lw1987}
D. Lewin,
\textit{Generalized musical intervals and transformations},
Yale University Press,
New Haven,
CT,
1987.
%
\bibitem{MAGMA2026}
MAGMA Developer Group,
\textit{MAGMA Computational Algebra System},
Version 2.29-10,
Sydney,
2026, \\
\texttt{http://magma.maths.usyd.edu.au}.
%
\bibitem{SA2015}
A. R. Santa Ana
\textit{Group theoretical approach on Riemannian uniform triadic transformations}, \\
Preprint,
2015,
Institute of Mathematics,
University of the Philippines,
Diliman.
%
\bibitem{OEIS2014}
N. J. A. Sloane,
\textit{The On-Line Encyclopedia of Integer Sequences} (OEIS), \\
The OEIS Foundation Inc.,
1964, \\
\texttt{http://oeis.org/}.
%
\bibitem{Wd2009}
C. Wood,
\textit{Abstracting tonality: Triads and uniform triadic transformations in an atonal context}, \\
Preprint,
2009,
Institute of Mathematics,
University of Chicago, \\
\texttt{http://www.math.uchicago.edu/\(\sim\)may/VIGRE/VIGRE2009/REUPapers/Wood.pdf}.
%
\bibitem{Zh2009}
A. Zhang,
\textit{The framework of music theory as represented with groups}, \\
Preprint,
2009,
Institute of Mathematics,
University of Washington. \\
\texttt{https://www.researchgate.net/publication/237773040}.
%
\end{thebibliography}
\end{document}